\documentclass{AIMS}
\usepackage{amsmath}
  \usepackage{paralist}
  \usepackage{graphics} %% add this and next lines if pictures should be in esp format
  \usepackage{epsfig} %For pictures: screened artwork should be set up with an 85 or 100 line screen
 \usepackage[colorlinks=true,breaklinks]{hyperref}
\hypersetup{urlcolor=blue, citecolor=red}
\def\publname{To appear in the proceedings of HYP2012}
\def\currentvolume{X}
 \def\currentissue{X}
  \def\currentyear{20XX}
   \def\currentmonth{XX}
    \def\ppages{X--XX}
     \def\DOI{}
\newtheorem{theorem}{Theorem}[section]

\newtheorem{proposition}{Proposition}

\theoremstyle{definition}

\newtheorem{remark}{Remark}

\newcommand{\ep}{\varepsilon}

\title[Semi-Lagrangian schemes for HJB equations]
      {Semi-Lagrangian schemes for linear and fully non-linear Hamilton-Jacobi-Bellman equations}

\author[Kristian Debrabant and Espen R.\ Jakobsen]{}


 \email{debrabant@imada.sdu.dk}
 \email{erj@math.ntnu.no}

\begin{document}
\maketitle

% Enter the first author's name and address:
\centerline{\scshape Kristian Debrabant }
\medskip
{\footnotesize
% please put the address of the first author
 \centerline{University of Southern Denmark, Department of Mathematics and Computer Science}
   \centerline{Campusvej 55}
   \centerline{5230 Odense M, Denmark}
} % Do not forget to end the {\footnotesize by the sign }

\medskip

\centerline{\scshape Espen Robstad Jakobsen}
\medskip
{\footnotesize
 % please put the address of the second  and third author
 \centerline{ Norwegian University of Science and Technology}
%   \centerline{Other lines}
   \centerline{NO--7491, Trondheim, Norway}
}

\bigskip

% The name of the associate editor will be entered by an editorial staff
% "Communicated by the associate editor name" is not needed for special issue.
 \centerline{(Communicated by the associate editor name)}

%The abstract of your paper
\begin{abstract}
We consider the numerical solution of Hamil\-ton-Jacobi-Bellman
equations arising in stochastic control theory. We introduce a class of
monotone approximation schemes relying on monotone
interpolation. These schemes converge under very weak assumptions,
including the case of arbitrary degenerate diffusions. Besides providing a
unifying framework that includes several
known first order accurate schemes, stability and convergence
results are given, along with two different robust error
estimates. Finally, the method is applied to a super-replication
problem from finance.
\end{abstract}

%The title of your section 1
\section{Introduction}\label{sec:local}
In this paper we consider the numerical solution of partial
differential equations of Hamil\-ton-Jacobi-Bellman type,
\begin{align}
\label{E}
u_t-\inf_{\alpha\in \ensuremath{\mathcal{A}}}\Big\{L^{\alpha} [u](t,x) +
c^{\alpha}(t,x) u +
f^{\alpha} (t,x)\Big\}&=0&&\text{in}\quad Q_T,\\
u(0,x)&=g(x)&&\text{in}\quad\ensuremath{\mathbb{R}}^N,
\label{IV}
\end{align}
where
\begin{align*}
&L^{\alpha} [u](t,x)= \mathrm{tr} [a^{\alpha}(t,x) D^2u(t,x)] +
b^{\alpha}(t,x) D
u(t,x),
\end{align*}
$Q_T:=(0,T]\times\ensuremath{\mathbb{R}}^N$, and $\ensuremath{\mathcal{A}}$ is a complete metric space.
The coefficients
$a^{{\alpha}}=\frac12\sigma^{\alpha}\sigma^{\alpha\,\top}$,
$b^{{\alpha}}$,
$c^{{\alpha}}$, $f^{{\alpha}}$ and the initial data $g$ take values
respectively in
$\ensuremath{\mathbb{S}}^N$, the space of $N \times N$ symmetric matrices, $\ensuremath{\mathbb{R}}^N$, $\ensuremath{\mathbb{R}}$,
$\ensuremath{\mathbb{R}}$, and $\ensuremath{\mathbb{R}}$. We will only assume that $a^{\alpha}$ is
positive semi-definite, thus the equation is allowed to degenerate and hence not have smooth solutions in general.
By solutions in this paper we will therefore always mean generalized solutions in the viscosity sense, see
e.\,g.\ \cite{CIL:UG,YZ}. Then the solution coincides with the value function of a finite horizon, optimal stochastic control problem \cite{YZ}.

To ensure comparison and well-posedness of \eqref{E}--\eqref{IV} in the class of bounded $x$-Lipschitz functions, we will use the following standard assumptions on its data:
\newcounter{KDERJcount:assumption}
\begin{enumerate}
\stepcounter{KDERJcount:assumption}
\renewcommand{\theenumi}{(A\arabic{KDERJcount:assumption})}
\renewcommand{\labelenumi}{\theenumi}
\item \label{A1}
 For any $\alpha\in \ensuremath{\mathcal{A}}$,
$a^{{\alpha}}=\frac12\sigma^{\alpha}\sigma^{\alpha\,\top}$ for some $N\times P$
matrix $\sigma^{{\alpha}}$. Moreover, there is a constant $K$ independent
of ${\alpha}$ such that
\[|g|_1+|\sigma^{{\alpha}}|_1+|b^{{\alpha}}|_1 + |c^{{\alpha}}|_1+|f^{{\alpha}}|_1 \leq
K,\]
where
 $|\phi|_1=\sup_{(t,x)\in Q_T}|\phi(x,t)|+\sup_{(x,t)\neq (y,s)}\frac{|\phi(x,t)-\phi(y,s)|}{|x-y|+|t-s|^{1/2}}$ is a space-time Lipschitz/H\"older-norm.
\end{enumerate}
The following result is standard.
\begin{proposition}
\label{WPi}
Assume that \ref{A1} holds. Then there exist a unique solution $u$ of
\eqref{E}--\eqref{IV} and a constant $C$ only depending on $T$ and $K$
from \ref{A1} such that
\[|u|_1\leq C.\]
Furthermore, if $u_1$ and $u_2$ are sub- and supersolutions of
\eqref{E} satisfying $u_1(0,\cdot)\leq u_2(0,\cdot)$,
then $u_1\leq u_2$.
\end{proposition}

\section{Semi-Lagrangian schemes}
Following \cite{debrabant13slsII} we propose a class of approximation  schemes
for \eqref{E}--\eqref{IV} which we call Semi-Lagrangian or SL
schemes. These schemes converge under very weak assumptions,
including the case of arbitrary degenerate diffusions. In particular, these
schemes are $L^\infty$-stable and convergent for problems involving
diffusion matrices that are not diagonally dominant.
This class includes (parabolic versions of) the
``control schemes'' of  Menaldi \cite{M} and Camilli and Falcone
\cite{CF:Appr} and some of the monotone schemes of Crandall and Lions
\cite{CL:MCM}. It also includes SL schemes for first
order Bellman equations \cite{CD,F} and some new versions as discussed in the following section.

The schemes are defined on a possibly unstructured family of grids
$\{G_{\ensuremath{\Delta t},\Delta x}\}$,
\[G=G_{\ensuremath{\Delta t},\Delta x}=\{(t_n,x_i)\}_{n\in\ensuremath{\mathbb{N}}_0,i\in\ensuremath{\mathbb{N}}}=\{t_n\}_{n\in\ensuremath{\mathbb{N}}_0}\times X_{\Delta x},\]
for $\ensuremath{\Delta t},\Delta x>0$. Here $0=t_0<t_1<\dots<t_n<t_{n+1}$ satisfy
\[\max_n\Delta t_n\leq\ensuremath{\Delta t} \qquad\text{where}\qquad \ensuremath{\Delta t}_n=t_{n}-t_{n-1} ,\]
and $X_{\Delta x}=\{x_i\}_{i\in\ensuremath{\mathbb{N}}}$ is the set of vertices or nodes for a
non-degenerate polyhedral subdivision of $\ensuremath{\mathbb{R}}^N$.

We consider the following general finite difference approximations of the differential operator $L^{\alpha}[\phi]$ in \eqref{E}:
\begin{align}
\label{BoAppr}
&L_k^{\alpha}[\phi](t,x):=
\sum_{i=1}^M\frac{\phi(t,x+y^{{\alpha},+}_{k,i}(t,x))-2\phi(t,x)+\phi(t,x+ y^{{\alpha},-}_{k,i}(t,x))}{2k^2},
\end{align}
for $k>0$ and some $M\geq1$.
For this approximation we will assume
{\makeatletter
\tagsleft@true
\makeatother
  \begin{align}\label{eq:Y1}\tag{Y1}
&\begin{cases}
&\displaystyle\sum_{i=1}^M [y^{{\alpha},+}_{k,i}+ y^{{\alpha},-}_{k,i}]=
  2k^2b^{\alpha}+\mathcal{O}(k^4),\\
&\displaystyle\sum_{i=1}^M [
y^{{\alpha},+}_{k,i}y^{{\alpha},+\,\top}_{k,i}+
  y^{{\alpha},-}_{k,i}y^{{\alpha},-\,\top}_{k,i}
  ]
= 2k^2
  \sigma^{\alpha}\sigma^{{\alpha}\,\top}+\mathcal{O}(k^4),\hspace{1.2cm}
  \\
&\displaystyle\sum_{i=1}^M [
y^{{\alpha},+}_{k,i,j_1}y^{{\alpha},+}_{k,i,j_2}y^{{\alpha},+}_{k,i,j_3}+
  y^{{\alpha},-}_{k,i,j_1}y^{{\alpha},-}_{k,i,j_2}y^{{\alpha},-}_{k,i,j_3}
]=\mathcal{O}(k^4),\\
&\displaystyle\sum_{i=1}^M [
y^{{\alpha},+}_{k,i,j_1}y^{{\alpha},+}_{k,i,j_2}y^{{\alpha},+}_{k,i,j_3}y^{{\alpha},+}_{k,i,j_4}+
  y^{{\alpha},-}_{k,i,j_1}y^{{\alpha},-}_{k,i,j_2}y^{{\alpha},-}_{k,i,j_3}y^{{\alpha},-}_{k,i,j_4}
]=\mathcal{O}(k^4),\\
\end{cases}\\
& \nonumber\quad\text{for all $j_1,j_2,j_3,j_4=1,2,\dots,N$ indicating
  components of the $y$-vectors.}
\end{align}}
Under assumption (\ref{eq:Y1}), a Taylor expansion shows that
$L_k^{\alpha}$ is a second order consistent approximation
satisfying
\begin{align}
\label{consistL}
|L^{\alpha}_k[\phi]-L^{\alpha}[\phi]|\leq C(|D\phi|_0+\dots+|D^4\phi|_0)k^2
\end{align}
for all smooth functions $\phi$, where $|\phi|_0=\sup_{(t,x)\in Q_T}|\phi(x,t)|$.

To relate this approximation to the spatial grid $X_{\Delta x}$, we replace $\phi$ by its interpolant $\mathcal{I}\phi$, yielding overall a
semi-discrete approximation of \eqref{E},
\begin{align*}
 U_t-\inf_{\alpha\in\mathcal A}\Big\{L_k^{\alpha}[\mathcal{I} U](t,x)+c^{\alpha}(t,x)U+f^{\alpha}(t,x)
\Big\}=0\quad
\text{in}\quad (0,T)\times X_{\Delta x}.
\end{align*}
We require the interpolation operator $\mathcal{I}$ to fulfill the following two conditions:
\begin{enumerate}
\renewcommand{\theenumi}{(I\arabic{enumi})}
\renewcommand{\labelenumi}{\theenumi}
\item \label{I1}
There are $K\geq0, r\in\ensuremath{\mathbb{N}}$ such
that for all smooth functions $\phi$
\[|(\mathcal{I}\phi)-\phi|_0\leq
K|D^r\phi|_0\Delta x^r.\]
\item\label{I2}There is a set of non-negative functions
  $\{w_{j}(x)\}_j$ such that
\[(\mathcal{I}\phi)(x)=\sum_j\phi(x_j) w_{j}(x),\]
and
\begin{gather*}
w_{j}(x)\geq0,\qquad w_{i}(x_j)=\delta_{ij}
\end{gather*}
for all $i,j\in \ensuremath{\mathbb{N}}$.
\end{enumerate}
\ref{I1} implies together with \eqref{consistL} that $L^{\alpha}_k[\mathcal{I}\phi]$ is a consistent approximation of $L^{\alpha}[\phi]$ if $\frac{\Delta x^r}{k^2}\to0$.  An interpolation  satisfying \ref{I2} is said to be {\em positive} and is {\em monotone} in the sense that $U\leq V$ implies that $\mathcal{I} U\leq\mathcal{I} V$. Typically $\mathcal{I}$ will be constant, linear, or multi-linear interpolation
(i.\,e.\ $r\leq 2$ in \ref{I1}), because higher order interpolation is not monotone in general.

The final scheme can now be found by discretizing in time using a parameter
$\theta\in[0,1]$,
\begin{align}
\label{FD}
&\delta_{\ensuremath{\Delta t}_n}U^n_i=\inf_{\alpha\in\mathcal{A}}\Big\{
L_k^{\alpha}[\mathcal{I}\bar{U}^{\theta,n}_\cdot]^{n-1+\theta}_i+c^{{\alpha},n-1+\theta}_i\bar{U}^{\theta,n}_i
+f^{{\alpha},n-1+\theta}_i\Big\}
\end{align}
in $G$, where $U^n_i=U(t_n,x_i)$,
  $f^{{\alpha},n-1+\theta}_i= f^{\alpha}(t_{n-1}+\theta\ensuremath{\Delta t}_n,x_i)$, $\dots$ for $(t_n,x_i)\in
  G$,
\[\delta_{\ensuremath{\Delta t}}\phi(t,x)=\frac{\phi(t,x)-\phi(t-\ensuremath{\Delta t},x)}{\ensuremath{\Delta t}},\qquad\text{and}
\qquad
\bar{\phi}^{\theta,n}_\cdot=(1-\theta)\phi^{n-1}_\cdot+\theta\phi^n_\cdot.\]
 As initial conditions we take
\begin{align}
&U_i^0=g(x_i)\quad\text{in}\quad X_{\Delta x}.\label{FD_BC}
\end{align}
For the choices $\theta=0,1$, and $1/2$ the time discretization corresponds to respectively explicit Euler, implicit Euler, and midpoint rule. For $\theta=1/2$, the full scheme can be seen as generalized Crank-Nicolson type discretization.

\section{Examples of approximations \texorpdfstring{$L^{\alpha}_k$}{}}
\label{exLalp}
\begin{enumerate}
\item The approximation of Falcone \cite{F} (see also \cite{CD}),
\[b^{{\alpha}}D\phi \approx \frac{\mathcal{I}\phi(x+hb^{\alpha})-\mathcal{I}\phi(x)}{h},\]
corresponds to our $L^{\alpha}_k$ if $k=\sqrt h$, $y_k^{{\alpha},\pm}=k^2b^{\alpha}$.

\item\label{Appr2} The approximation of Crandall-Lions \cite{CL:MCM},
\[\frac12\mathrm{tr}[\sigma^{{\alpha}}\sigma^{{\alpha}\,\top}D^2\phi]\approx \sum_{j=1}^P\frac{\mathcal{I}\phi(x+k\sigma_j^{\alpha})-2\mathcal{I}\phi(x)+\mathcal{I}\phi(x-k\sigma_j^{\alpha})}{2k^2},\]
corresponds to our $L^{\alpha}_k$ if $y_{k,j}^{{\alpha},\pm}=
\pm k\sigma_j^{\alpha}$ and $M=P$.

\item The corrected version of the approximation of Camilli-Falcone
\cite{CF:Appr} (see also \cite{M}),
\begin{align*}
&\frac12\mathrm{tr}[\sigma^{{\alpha}}\sigma^{{\alpha}\,\top}D^2\phi]+b^{{\alpha}}D\phi\\
& \approx
\sum_{j=1}^P\frac{\mathcal{I}\phi(x+\sqrt{h}\sigma_j^{\alpha}+\frac{h}{P}b^{\alpha})-2\mathcal{I}\phi(x)+\mathcal{I}\phi(x-\sqrt{h}\sigma_j^{\alpha}+\frac{h}{P}b^{\alpha})}{2h},
\end{align*}
corresponds to our $L^{\alpha}_k$ if $k=\sqrt h$, $y_{k,j}^{{\alpha},\pm}=
\pm k\sigma_j^{\alpha}+\frac{k^2}{P}b^{\alpha}$  and $M=P$.

\item\label{Appr4} The new approximation obtained by combining approximations 1 and 2,
\begin{align*}
&\frac12\mathrm{tr}[\sigma^{{\alpha}}\sigma^{{\alpha}\,\top}D^2\phi]+b^{{\alpha}}D\phi\\
&\approx \frac{\mathcal{I}\phi(x+k^2b^{\alpha})-\mathcal{I}\phi(x)}{k^2}+
\sum_{j=1}^P\frac{\mathcal{I}\phi(x+k\sigma_j^{\alpha})-2\mathcal{I}\phi(x)
  +\mathcal{I}\phi(x-k\sigma_j^{\alpha})}{2k^2},
\end{align*}
corresponds to our $L^{\alpha}_k$ if $y_{k,j}^{{\alpha},\pm}= \pm
k\sigma_j^{\alpha}$ for $j\leq P$, $y_{k,P+1}^{{\alpha},\pm}= k^2 b^{\alpha}$ and $M=P+1$.
\item\label{Appr5} Yet another new approximation,
\begin{multline*}
\frac12\mathrm{tr}[\sigma^{{\alpha}}\sigma^{{\alpha}\,\top}D^2\phi]+b^{{\alpha}}D\phi
\approx
\sum_{j=1}^{P-1}\frac{\mathcal{I}\phi(x+k\sigma_j^{\alpha})-2\mathcal{I}\phi(x)+\mathcal{I}\phi(x-k\sigma_j^{\alpha})}{2k^2}\\
+\frac{\mathcal{I}\phi(x+k\sigma_P^{\alpha}+k^2b^{\alpha})-2\mathcal{I}\phi(x)+\mathcal{I}\phi(x-k\sigma_P^{\alpha}+k^2b^{\alpha})}{2k^2},
\end{multline*}
corresponds to our $L^{\alpha}_k$ if $y_{k,j}^{{\alpha},\pm}= \pm
k\sigma_j^{\alpha}$ for $j<P$, $y_{k,P}^{{\alpha},\pm}= \pm k\sigma_P^{\alpha}+k^2 b^{\alpha}$ and $M=P$.
\end{enumerate}
When $\sigma^{\alpha}$ does not depend on ${\alpha}$ but $b^{\alpha}$ does, approximations 4 and 5 are much more efficient than approximation 3.

\section{Linear interpolation SL scheme (LISL)}
To keep the scheme \eqref{FD} monotone, linear or
multi-linear interpolation is the most accurate interpolation
one can use in general. In this typical case we call the full scheme
\eqref{FD}--\eqref{FD_BC} the LISL scheme. In the following, we denote
by $c^{\alpha,+}$ the positive part of $c^{\alpha}$. Then we have the following result by \cite{debrabant13slsII}:

\begin{theorem}\label{cor:LISL}
Assume that \ref{A1}, \ref{I1}, \ref{I2}, and (\ref{eq:Y1}) hold.
\smallskip

\noindent (a) The LISL scheme is monotone if the following CFL
conditions hold:
\begin{align}
\label{CFL}
 (1-\theta)\ensuremath{\Delta t}\Big[\frac
  {M}{k^2}-c^{{\alpha},n-1+\theta}_{i}\Big]\leq
  1\ \ \text{and}\ \ \theta\ensuremath{\Delta t}\,
  c^{{\alpha},n-1+\theta}_{i}\leq1 \ \text{for all}\ {\alpha},n,i.
\end{align}

\noindent (b) The truncation error of the LISL scheme is
$O(|1-2\theta|\ensuremath{\Delta t} +\ensuremath{\Delta t}^2+ k^2 +\frac{\Delta x^2}{k^2} )$;
it is first order accurate for
$k=O(\ensuremath{\Delta x}^{1/2})$, $\ensuremath{\Delta t}=O(\Delta x)$ ($\ensuremath{\Delta t}=O(\Delta x^{1/2})$ if $\theta=\frac12$).
\smallskip

\noindent (c)  If
$2\theta\ensuremath{\Delta t}\sup_{\alpha}|c^{\alpha,+}|_0\leq1$ and
\eqref{CFL} holds, then there exists a unique bounded and
$L^\infty$-stable solution $U$ of the LISL scheme converging
uniformly to the solution $u$ of \eqref{E}--\eqref{IV} as
$\ensuremath{\Delta t},k,\frac{\Delta x}k\rightarrow 0$.
\end{theorem}

From this result it follows that the scheme is at most {\em first
  order accurate}, has {\em wide and increasing stencil} and a {\em
  good CFL condition}.  From the truncation error and the definition
of $L_k^{\alpha}$ the stencil is wide since the scheme is
consistent only if $\Delta x/k\rightarrow0$ as $\Delta x\rightarrow0$ and has stencil length
proportional to
  \[l:=\frac{\underset{t,x,{\alpha},i}{\max}\{ |y^{{\alpha},-}_{k,i}|,
    |y^{{\alpha},+}_{k,i}|\}}{\Delta x}\sim \frac
  k{\Delta x}\rightarrow\infty\quad\text{as}\quad\Delta x\rightarrow0.\] Here we have used that
  if (\ref{eq:Y1}) holds and $\sigma\not\equiv 0$, then typically
  $y^{{\alpha},\pm}_{k,i}\sim k$.  Note that if
  $k=\Delta x^{1/2}$, then $l\sim \Delta x^{-1/2}.$ Finally, in the case
  $\theta\neq1$ the CFL condition for \eqref{FD} is $\ensuremath{\Delta t}\leq
  Ck^2\sim\Delta x$ when $k=O(\Delta x^{1/2})$, and it is much less restrictive
  than the usual parabolic CFL condition, $\ensuremath{\Delta t}=O(\Delta x^2)$.

\begin{remark}
The LISL scheme is consistent and monotone for arbitrary degenerating diffusions, without requiring that $a^\alpha$
is diagonally dominant or similar conditions. In comparison to other schemes applicable in this situation,
 like the ones of Bonnans-Zidani \cite{BOZ}, it is much easier to analyze and to implement and faster in the sense that
 the computational cost for approximating the diffusion matrix is for fixed $x,t,\alpha$ independent of the stencil size.
\end{remark}
\section{The error estimate of \texorpdfstring{\cite{debrabant13slsII}}{Debrabant and Jakobsen (2013)}}
\label{Sec:ErrBnd}
To simplify the presentation, in the following we restrict to
a uniform time-grid, $G=\ensuremath{\Delta t}\,\{0,1,\dots,N_T\}\times X_{\Delta x}$. Let
$Q_{\ensuremath{\Delta t},T}:=\ensuremath{\Delta t}\,\{0,1,\dots,N_T\}\times\ensuremath{\mathbb{R}}^N$.
To apply the regularization method of Krylov \cite{Kr:00} we need a regularity
and continuous dependence result for the scheme that relies on the following additional (covariance-type)
assumptions: Whenever two sets of data $\sigma, b$ and $\tilde{\sigma}, \tilde{b}$ are
given, the corresponding approximations $L_k^\alpha, y^{\alpha,\pm}_{k,i}$
and $\tilde{L}_k^\alpha, \tilde{y}^{\alpha,\pm}_{k,i}$ in \eqref{BoAppr} satisfy
{\makeatletter
\tagsleft@true
\makeatother
\begin{align}\label{eq:Y3}
\tag{Y2}
\begin{cases}
  &\displaystyle\sum_{i=1}^M [y^{\alpha,+}_{k,i}+
  y^{\alpha,-}_{k,i}]-[\tilde{y}^{\alpha,+}_{k,i}+ \tilde{y}^{\alpha,-}_{k,i}]\leq
  2k^2(b^\alpha-\tilde b^\alpha),\\
  &\displaystyle\sum_{i=1}^M [
 y^{\alpha,+}_{k,i}y^{\alpha,+\,\top}_{k,i}+ y^{\alpha,-}_{k,i}y^{\alpha,-\,\top}_{k,i}
 ]+ [
  \tilde{y}^{\alpha,+}_{k,i}\tilde{y}^{\alpha,+\,\top}_{k,i}+\tilde{y}^{\alpha,-}_{k,i}\tilde{y}^{\alpha,-\,\top}_{k,i}
  ]\hspace{1cm}\\
  &\displaystyle\qquad -[
  y^{\alpha,+}_{k,i}\tilde{y}^{\alpha,+\,\top}_{k,i}+\tilde{y}^{\alpha,+}_{k,i}y^{\alpha,+\,\top}_{k,i}+y^{\alpha,-}_{k,i}\tilde{y}^{\alpha,-\,\top}_{k,i}+\tilde{y}^{\alpha,-}_{k,i}y^{\alpha,-\,\top}_{k,i}
] \\
  &\displaystyle \leq
  2k^2(\sigma^\alpha-\tilde{\sigma}^\alpha)(\sigma^\alpha-\tilde{\sigma}^\alpha)^\top +
  2k^4(b^\alpha-\tilde{b}^\alpha)(b^\alpha-\tilde{b}^\alpha)^\top,
\end{cases}
\end{align}}
when $\sigma,b,y_k^\pm$ are evaluated at $(t,x)$ and $\tilde\sigma,\tilde
b,\tilde y_k^\pm$ are evaluated at $(t,y)$ for all $t,x,y$.

Then one can prove the following error estimate \cite{debrabant13slsII}:
\begin{theorem}[Error Bound I]
Assume that \ref{A1}, \ref{I1}, \ref{I2},
(\ref{eq:Y1}), and
 (\ref{eq:Y3}) hold, and that $\ensuremath{\Delta t},\Delta x>0$, $k\in(0,1)$ satisfy the CFL conditions \eqref{CFL}.  If $u$ solves
\eqref{E}--\eqref{IV} and $U$ solves \eqref{FD}--\eqref{FD_BC}, then
there is $c_0>0$ such that for any $\ensuremath{\Delta t}\in(0,c_0)$
\[|u-U|\leq C(|1-2\theta|\ensuremath{\Delta t}^{1/4}+\ensuremath{\Delta t}^{1/3}+k^{1/2}+\frac{\Delta x}{k^2})\quad\text{in}\quad G.\]
\end{theorem}

This error bound holds also for unstructured grids. For  more regular solutions
it is possible to obtain better error estimates, but general
and optimal results are not available.
The best estimate in our case is $O(\Delta x^{1/5})$ which is achieved when
$k=O(\Delta x^{2/5})$ and $\ensuremath{\Delta t}=O(k^2)$. Note that the CFL conditions
\eqref{CFL} already imply that $\ensuremath{\Delta t}=O(k^2)$ if
$\theta<1$.  Also note that the above bound does not show convergence
when $k$ is optimal for the LISL scheme ($k=O(\Delta x^{1/2})$).

\section{A new error estimate}
In the above error estimate, the lower estimate on $u-U$ follows if you
can prove regularity and continuous dependence results for the
solution of the equation only. The proof of the upper estimate is
symmetric and requires such results for the numerical
solution. However, it is possible to avoid using such properties of
the numerical solution by a clever approximation argument, see
e.\,g.\ \cite{BJ:Rate3}. This allows for error estimates that show convergence for
any $k$ such that the scheme is consistent. We need an extra assumption on the coefficients:
\smallskip

\begin{enumerate}
\stepcounter{KDERJcount:assumption}
\renewcommand{\theenumi}{(A\arabic{KDERJcount:assumption})}
\renewcommand{\labelenumi}{\theenumi}
\item \label{A2}
The coefficients $\sigma^\alpha$ , $b^\alpha$,
  $c^\alpha$, $f^\alpha$ are continuous in $\alpha$ for all $x, t$.
\end{enumerate}

\begin{theorem}[Error Bound II]
Assume that \ref{A1}, \ref{A2}, \ref{I1} with $r=2$ ($\sim$linear
interpolation), \ref{I2}, and
(\ref{eq:Y1}) hold, and that $\ensuremath{\Delta t},\Delta x>0$, $k\in(0,1)$ satisfy the CFL conditions \eqref{CFL}. If $u$ solves
\eqref{E}--\eqref{IV} and $U$ solves \eqref{FD}--\eqref{FD_BC}, then
there is $c_0>0$ such that for any $\ensuremath{\Delta t}\in(0,c_0)$
\begin{align*}
u-U &\geq
C\Big(|1-2\theta|\ensuremath{\Delta t}^{1/4}+\ensuremath{\Delta t}^{1/3}+k^{1/2}+\frac{\Delta x}{k}\Big)\quad\text{in}\quad
G,\\
u-U &\leq
C\Big(|1-2\theta|\ensuremath{\Delta t}^{1/10}+\ensuremath{\Delta t}^{1/8}+k^{1/5}+\big(\frac{\Delta x}{k}\big)^{1/2}\Big)\quad\text{in}\quad G.
\end{align*}
\end{theorem}
With optimal $k$ for the LISL
scheme, $\ensuremath{\Delta t}=O(k^2)$ and $k=O(\Delta x^{1/2})$, we find that
$u-U=O(\Delta x^{1/10})$.

\begin{proof}
By a direct computation the local truncation error of
the method is bounded by
\begin{align*}
\frac{|1-2\theta|}2|\phi_{tt}|_0\ensuremath{\Delta t}
+C\Big(&\ensuremath{\Delta t}^2\left(
|\phi_{tt}|_0
+
|\phi_{ttt}|_0+
|D\phi_{tt}|_0+|
D^2\phi_{tt}|_0
\right)
\\
&+ |D^2\phi|_0
 \frac{\Delta x^2}{k^2}+(|D\phi|_0+\dots+|D^4\phi|_0)k^2
 \Big)
\end{align*}
for smooth $\phi$ (cf.\ Lemma 4.1 in
\cite{debrabant13slsII}). Moreover if also
$\partial_t^{k_1}D_x^{k_2}\phi=O(\ep^{1-2k_1-k_2})$ for any $k_1,k_2\in\ensuremath{\mathbb{N}}_0$, then the truncation
error is of order
\[(1-2\theta)\ensuremath{\Delta t}
\ep^{-3}+\ensuremath{\Delta t}^2\ep^{-5}+k^2\ep^{-3}+\frac{\Delta x^2}{k^2}\ep^{-1}=:E(\ep).\]
Since the scheme is monotone (under the CFL condition) and condition
\ref{A1} holds, it now follows from Theorem 3.1 in \cite{BJ:Rate3} that
\[C\inf_{\ep>0}\big(\ep+E(\ep)\big)\leq u- U\leq
C\inf_{\ep>0}\big(\ep^{1/3}+E(\ep)\big),\]
and we complete the proof optimizing over $\ep$ (as e.\,g.\ in \cite{BJ:Rate3,debrabant13slsII}).
\end{proof}

\section{Convergence test for a super-replication problem}\label{subsec:Bruder}
We consider a test problem from \cite{BBMZ} which was used to test
convergence rates for numerical approximations of a super-replication
problem from finance. The corresponding PDE is
\begin{equation}\label{eq:Superrepeq}
\inf_{\alpha_1^2+\alpha_2^2=1}
\left\{
\alpha_1^2u_t(t,x)-\frac12\mathrm{tr}\left(\sigma^\alpha(t,x)\sigma^{\alpha\,\top}(t,x) D^2u(t,x)\right)
\right\}
=f(t,x)
\end{equation}
with $0\leq x_1,x_2\leq 3$, $\sigma^\alpha(t,x)=\begin{pmatrix}
\alpha_1x_1\sqrt{x_2}\\
\alpha_2\eta(x_2)
\end{pmatrix}$ and $\eta(x)=x(3-x)$. We take
$u(t,x)=1+t^2-e^{-x_1^2-x_2^2}$ as exact solution as in \cite{BBMZ},
and then $f$ is forced to be
\begin{multline*}
f(t,x)=
\frac12
\left(
u_t-\frac12x_1^2x_2u_{x_1x_1}-\frac12x_2^2(3-x_2)^2u_{x_2x_2}
\right.\\
\left.
-\sqrt{\left(-u_t+\frac12x_1^2x_2u_{x_1x_1}-\frac12x_2^2(3-x_2)^2u_{x_2x_2}\right)^2+
\left(x_1\sqrt{x_2}^3(3-x_2)u_{x_1x_2}\right)^2
}
\right).
\end{multline*}
In \cite{BBMZ} $\eta(x)=x$, while we take $\eta(x)=x(3-x)$
to prevent the LISL scheme from overstepping the boundaries.
Note that changing $\eta$ does {\em not} change the solutions as long
as $\eta>0$ in the interior of the domain, see \cite{BBMZ}, and hence
the above equation is equivalent to the equation used in \cite{BBMZ}.
The initial values and Dirichlet boundary values at $x_1=0$ and $x_2=0$ are
taken from the exact solution.  As
in \cite{BBMZ}, at $x=3$ and $y=3$ homogeneous Neumann boundary
conditions are implemented.  To approximate the values of
$\alpha_1,\alpha_2$, the Howard algorithm is used (see \cite{BBMZ}),
which requires an implicit time discretization, so we choose
$\theta=1$.
We choose $k=\sqrt{\Delta x}$ and a regular triangular grid. The numbers of time steps are chosen as $\frac1{\ensuremath{\Delta x}}$.

The results at $t=1$ are given in Table \ref{Tab:Ex7bundbc}.  The
numerical order of convergence is approximately one.

\begin{table}[ht]
\begin{tabular}[t]{c|c|c}
$\Delta x$&$|u-U|_0$&rate\\\hline
1.50e-1&2.01e-1&\\
7.50e-2&9.49e-2&1.08\\
3.75e-2&4.29e-2&1.15\\
1.87e-2&1.94e-2&1.15
\end{tabular}
}
\caption{Results for the convergence test for the super-replication problem at $t=1$}{\label{Tab:Ex7bundbc}
\end{table}
\begin{remark}
Equation \eqref{eq:Superrepeq} can not be written in a form
\eqref{E} satisfying the assumptions of this paper, so the results of
this paper do not apply to this problem. However, it seems possible to
extend them to cover this problem using comparison results from \cite{BBMZ} along with
$L^\infty$-bounds on the numerical solution that follow
from the maximum principle.
\end{remark}
\section{A super-replication problem}\label{subsec:superrep}
\begin{figure}
\includegraphics{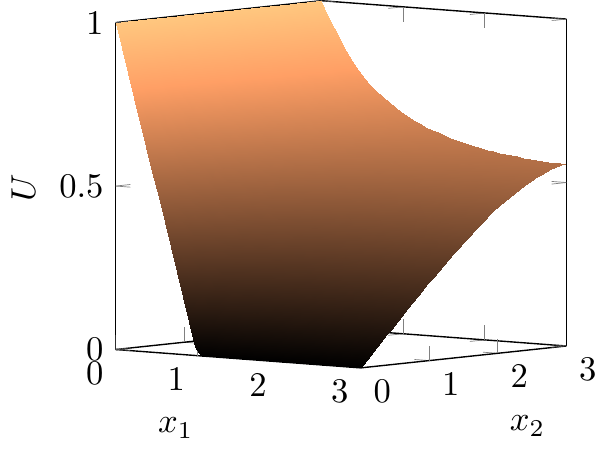}
\caption{Numerical solution of super-replication problem at $t=1$}\label{Fig:Ex6}
\end{figure}

We apply our method to solve a problem from finance, the
super-replication problem under gamma constraints considered in
\cite{BBMZ}. It consists of solving equation \eqref{eq:Superrepeq}
with $f\equiv0$, Neumann boundary conditions
and $\sigma^\alpha$ as in Subsection \ref{subsec:Bruder}, and
initial and Dirichlet conditions given by
\begin{align*}
u(t,x)=&\max(0,1-x_1),\quad t=0\quad or \quad x_1=0\quad or \quad x_2=0.
\end{align*}
The solution obtained with the LISL scheme is given in Figure
\ref{Fig:Ex6} and coincides with the solution found in \cite{BBMZ}. It
gives the price of a put option of strike and
maturity 1, and $x_1$ and $x_2$ are respectively the price of the
underlying and the price of the forward variance swap on the
underlying.
\newlength{\figwidth}
\setlength{\figwidth}{0.485 \textwidth}

\medskip
% The data information below will be filled by AIMS editorial staff
Received xxxx 20xx; revised xxxx 20xx.
\medskip

\end{document}